\documentclass{article}
\RequirePackage{etex} %
\usepackage[letterpaper,margin=1in]{geometry}
\usepackage{amsmath,amssymb,amsthm,amsfonts,latexsym,bbm,xspace,graphicx,float,mathtools,mathdots}
\usepackage{braket,caption,subcaption,ellipsis,xcolor,textcomp}
\usepackage{combelow} %
\usepackage{booktabs}
\usepackage{hyperref}
\usepackage[nameinlink]{cleveref}
\crefname{ineq}{inequality}{inequalities}
\creflabelformat{ineq}{#2{\upshape(#1)}#3}

\usepackage{enumitem}

\newtheorem{theorem}{Theorem}

\usepackage{chngcntr}
\counterwithin{theorem}{section}

\newtheorem*{namedtheorem}{\theoremname}
\newcommand{\theoremname}{testing}

\newtheorem{problem}[theorem]{Problem}

\newtheorem{proposition}[theorem]{Proposition}

\newtheorem{assumption}[theorem]{Assumption}

\theoremstyle{definition}
\newtheorem{definition}[theorem]{Definition}
\newtheorem{remark}[theorem]{Remark}

\newcommand{\bone}{\boldsymbol{1}}

\newcommand{\ignore}[1]{}

\usepackage{autonum}

\usepackage{natbib}
\setcitestyle{authoryear,open={(},close={)}} %

\DeclareMathOperator*{\argmax}{argmax}

\newcommand{\iid}{\text{i.i.d.}\xspace}

\newcommand{\mc}[1]{\mathcal{#1}}

\newcommand{\lb}{\left[}
\newcommand{\rb}{\right]}
\newcommand{\lp}{\left(}
\newcommand{\rp}{\right)}

\newcommand{\defined}{\coloneqq}

\newcommand{\cusum}{\text{CuSum}\xspace}

\newcommand{\dkl}{d_{\text{KL}}}

\newcommand{\blue}[1]{#1}

\newcommand{\bcsdetector}{\texttt{BCS-Detector}\xspace}
\newcommand{\rfcsdetector}{\texttt{repeated-FCS-Detector}\xspace}
\newcommand{\filtration}{\mathcal{F}}

\usepackage{pgfplots}
\pgfplotsset{compat=newest}
\pgfplotsset{scaled y ticks=false}
\usepgfplotslibrary{groupplots}
\usepgfplotslibrary{dateplot}
\tikzstyle{every node}=[font=\small]
\pgfplotsset{
    yticklabel style={/pgf/number format/fixed},  
}

\pgfplotsset{compat=1.11,
 /pgfplots/ybar legend/.style={
 /pgfplots/legend image code/.code={
 \draw[##1,/tikz/.cd,yshift=-0.25em]
 (0cm,0cm) rectangle (3pt,0.8em);},
 },
}

\title{Reducing sequential change detection to sequential estimation}
\usepackage{authblk}
\author[1]{Shubhanshu Shekhar}
\author[1, 2]{Aaditya Ramdas}
\affil[1]{Department of Statistics and Data Science, Carnegie Mellon University}
\affil[2]{Machine Learning Department, Carnegie Mellon University}
\affil[ ]{\texttt{\{shubhan2,aramdas\}@andrew.cmu.edu}}

\begin{document}
\date{}
\maketitle

\begin{abstract}
    We consider the problem of sequential change detection, where the goal is to design a scheme for detecting any changes in a parameter or functional $\theta$ of the data stream distribution that has small detection delay, but guarantees control on the frequency of false alarms in the absence of changes. In this paper, we describe a simple reduction from sequential change detection to sequential estimation using confidence sequences: we begin a new $(1-\alpha)$-confidence sequence at each time step, and proclaim a change when the intersection of all active confidence sequences becomes empty. We prove that the average run length is at least $1/\alpha$, resulting in a change detection scheme with minimal structural assumptions~(thus allowing for possibly dependent observations, and nonparametric distribution classes), but strong guarantees. 
    We also describe an interesting parallel with Lorden's reduction from change detection to sequential testing and connections to the recent ``e-detector'' framework. 
\end{abstract}

\section{Introduction}
\label{sec:introduction}
    We consider the following problem of \emph{sequential change detection} (SCD): for a general space $\mc{X}$, given a stream of $\mc{X}$-valued observations $X_1, X_2, \ldots$, our goal is to design a method to detect any changes in a prespecified parameter~or functional $\theta$ (possibly infinite dimensional) associated with the source generating this stream. Let $\mc{P}$ denote a class of probability distributions on the infinite product space $\Omega = \mc{X}^{\infty}$, and let $\theta: \mc{P} \to \Theta$, denote a mapping from probability distributions to some (possibly infinite dimensional) parameter space $\Theta$. The data distribution satisfies the following for some $T \in \mathbb{N}\cup\{\infty\}$: 
    \begin{itemize}
        \item For $n \leq T$, the observations are the first $T$ elements of a trajectory of $P_0 \in \mc{P}$, with $\theta(P_0) = \theta_0$. 
        \item For $n > T$, the observations are drawn from another distribution $P_1 \in \mc{P}$, such that $\theta(P_1) = \theta_1 \neq \theta_0$.  
    \end{itemize} 
    In particular, $P_0, P_1$ are distributions over an infinite sequence of observations, and the data is not assumed to be \iid (independent and identically distributed), meaning that we do not assume that $P_0,P_1$ take the form $p^\infty$ for some distribution $p$ over $\mathcal X$. Further, $P_1$ (and hence $\theta_1$) is allowed to depend on $X_1,\dots,X_T$.
    
    Under the above model, our goal is to design a data-driven method for detecting any change in the value of $\theta$, without any knowledge of $T,P_0,P_1$. In technical terms, our objective is to define a stopping time $\tau$, at which we stop collecting more data, and declare a changepoint has previously occurred. 

    There are two problem settings to consider: non-partitioned and partitioned. By default, define the filtration $\filtration\equiv (\mc{F}_n)_{n \geq 0}$ as given by $\mc{F}_n \defined \sigma(X_1, \ldots, X_n)$ and $\mc{F}_0 = \{\emptyset, \Omega\}$.
    
    \begin{problem}[Non-partitioned SCD]
        \label{problem:scd} For some unknown triple $(T, P_0, P_1)$, suppose $X_1, X_2, \ldots$ denote a data stream, such that $(X_n)_{n \leq T}$ are drawn from $P_0$, and $(X_n)_{n>T}$ are drawn from $P_1$, such that $\theta(P_1) = \theta_1 \neq \theta_0 = \theta(P_0)$. We must define a stopping time $\tau$, adapted to $\filtration$ satisfying: 
        \begin{itemize}
            \item If $T=\infty$, and there is no changepoint, we require that $\mathbb{E}_{\infty}[\tau] \geq 1/\alpha$, for a prespecified $\alpha \in (0,1)$. The term $\mathbb{E}_{\infty}[\tau]$ is called the average run length~(ARL), and represents the frequency of false alarms. 
            \item If $T<\infty$, and there is a changepoint at which the distribution changes from $P_0$ to $P_1$, we desire the detection delay, $\mathbb{E}_{T}[(\tau-T)^+]$, to be as small as possible. 
        \end{itemize}
    \end{problem}
    
    The non-partitioned SCD problem stated above, does not require any pre-specified partitioning of $\mc{P}$ into pre- and post change distribution classes. That is, it assumes that the data generating distribution changes from one unknown distribution~($P_0$) in $\mc{P}$ to another distribution~($P_1$). This is in contrast to another formulation of the SCD problem, where it is assumed that we know some partition of $\mc{P}$ into $\mc{P}_0,\mc{P}_1$ such that $P_0,P_1$ are known to respectively lie in $\mc{P}_0,\mc{P}_1$ (in other words, we are not interested in changes within $\mc{P}_0$, only changes from $\mc{P}_0$ to $\mc{P}_1$).  This additional knowledge about the two partitioned distribution classes, that is not available in the first formulation, is often critical in designing optimal SCD schemes, especially in parametric problems. The strategy we develop in this paper is  also applicable to an intermediate variant of the SCD problem, that only assumes the knowledge of the pre-change parameter class, $\Theta_0$.
    \begin{problem}[Partitioned SCD]
        \label{problem:scd-separated}
        For some unknown triple $(T, P_0, P_1)$, suppose $X_1, X_2, \ldots$ denote a stream of $\mc{X}$-valued observations, satisfying the following assumptions:
        The observations $(X_n)_{n \leq T}$ are drawn according to a process $P_0 \in \mc{P}$ with parameter $\theta_0 \in \Theta_0 \subset \Theta$, and this parameter set $\Theta_0$ is assumed to be known. 
        The observations $(X_n)_{n > T}$ are drawn from $P_1$ with parameter  $\theta_1$, such that $\theta_1 \not \in \Theta_0$. The parameter $\theta_1$ is not assumed to be known. 
        Our objective is to design a stopping time $\tau$, satisfying the two bulleted criteria stated in~\Cref{problem:scd}. 
    \end{problem}

    Sequential changepoint detection is a very well-studied problem in the sequential analysis literature, going back to the early works by~\citet{shewhart1925application, shewhart1930economic, page1954continuous, shiryaev1963optimum}. These initial papers developed computationally efficient likelihood-based schemes for known pre- and post-change distributions, which have since then been extended to more general cases of composite, but parametric, pre- and post-change distribution classes. We refer the reader to the book by~\citet{tartakovsky2014sequential} for a detailed discussion of the parametric SCD problem. Unlike the parametric case, there exist very few general principles~(analogous to the likelihood and generalized likelihood based schemes) for designing SCD methods with nonparametric distribution classes. Two recent exceptions to this trend include the paper on e-detectors by~\citet{shin2022detectors} for~\Cref{problem:scd-separated}~(partitioned),  and the \bcsdetector scheme proposed by~\citet{shekhar23changepoint} for~\Cref{problem:scd}~(non-partitioned), which is the subject of this paper. 

    The \bcsdetector scheme of~\citet{shekhar23changepoint} is also a reduction from changepoint detection to estimation, but in this paper we propose a different, and even simpler reduction. To elaborate, \bcsdetector uses a single confidence sequence~(CS) in the forward direction, but with every new observation, it constructs a new CS in the backward direction~(the so-called ``backward CS'' or BCS, constructed using observations with their time indices reversed). The scheme stops and declares a detection, as soon as the intersection of the CS and the all active BCSs  becomes empty.  Since it is critical to our simplified scheme as well, we recall the definition of CSs below. 
    \begin{definition}[Confidence Sequences]
        \label{def:conf-seqs}
        Given observations $X_1, X_2, \ldots$ drawn from a distribution $P$ with associated parameter $\theta \equiv \theta(P)$, a level-$(1-\alpha)$ CS for  $\theta$, is a sequence of sets $(C_n)_{n \geq 1}$, satisfying: 
        \begin{itemize}
            \item For every $n \geq 1$, the set $C_n \subset \Theta$ is $\mc{F}_n = \sigma(X_1, \ldots, X_n)$-measurable. In words, the set $C_n$ can be constructed using the information contained in the first $n$ observations. 
            \item The sets satisfy a uniform coverage guarantee: $\mathbb{P}\lp \forall n \in \mathbb{N}: \theta(P)  \in C_n \rp \geq 1-\alpha$. Equivalently, for any $\mc{F}$-stopping time $\tau$, $\mathbb{P}\lp  \theta(P)  \in C_\tau \rp \geq 1-\alpha$.
        \end{itemize}
    \end{definition}

    \begin{remark}\label{rem:nested}
                Due to the uniform coverage guarantee, if $(C_n)$ is a CS, then so is $(\cap_{m \leq n} C_m)$. Thus, we can assume without loss of generality that the sets involved in a CS are nested; that is $C_n \subset C_{n'}$ for all $n'<n$.
    \end{remark}
    
    The \bcsdetector scheme of~\citet{shekhar23changepoint} satisfies several favorable properties: it can be instantiated for a large class of parametric and nonparametric problems, it provides non-asymptotic control over the ARL, and has strong guarantees over the detection delay. However, a closer look at their scheme reveals that it implicitly makes a ``bidirectional'' assumption about the data generating process: at any $n \geq 1$, the \bcsdetector assumes the ability to construct a CS  in the forward direction~(based on $X_1, \ldots, X_n$), as well as in the backward direction~(using $Y_1 = X_n, Y_2=X_{n-1},\ldots, Y_n = X_1$). Most methods for constructing CSs proceed by designing martingales or supermartingales adapted to the natural (forward) filtration of the observations. Hence, the \bcsdetector implicitly involves constructing martingales~(or supermartingales) in both, the forward and reverse directions; this in turn is typically only possible if the observations are independent. This restriction limits the applicability of the \bcsdetector scheme, and this paper's simplified scheme addresses this weakness.  
    
    \paragraph{Contributions.} Our main contribution in this paper is a new SCD scheme, that we refer to as the \rfcsdetector. This scheme proceeds by constructing a new forward CS with each observation~(\Cref{def:repeated-fcs-detector}), and stops as soon the intersection of all the active CSs becomes empty. Since it relies only on forward CSs, it eliminates the aforementioned weakness of the \bcsdetector of~\citet{shekhar23changepoint}. Further, as we show in~\Cref{prop:arl-control}, our scheme achieves a tighter bound~(by a factor $2$)  on the ARL, as compared to the \bcsdetector, while matching its detection delay guarantees. Finally, we note in~\Cref{sec:lorden}, that our reduction from change detection to constructing CSs  is 
    significantly generalizes a famous reduction by \cite{lorden1971procedures} from parametric change detection to one-sided sequential tests.

\section{Our proposed reduction}
\label{sec:proposed-scheme}
    We now describe our new SCD scheme that proceeds by starting a new CS in the forward direction with each new observation. 
    \begin{definition}[\rfcsdetector]
        \label{def:repeated-fcs-detector} Suppose we are given a stream of observations $X_1, X_2, \ldots$, and a functional $\theta$ associated with the source. For~\Cref{problem:scd}~(non-partitioned), 
         define the  $C^{(0)}_n = \Theta$ for all $n \geq 1$, while for~\Cref{problem:scd-separated}~(partitioned) set $C^{(0)}_n = \Theta_0$, for all $n \geq 1$.   Proceed as follows for $n=1, 2, \ldots:$
        \begin{itemize}
            \item Observe the next data-point $X_n$. 
            \item Initialize a new level-$(1-\alpha)$ forward CS, $C^{(n)} \equiv (C^{(n)}_t)_{t \geq n}$, which will be formed using data $X_n, X_{n+1},\dots$. 
            \item Update all  CSs initialized so far using $X_n$, meaning that we form the sets $\{C_n^{(m)}: 1 \leq m \leq n\}$.
            \item If the intersection of all initialized CSs becomes empty, meaning $\cap_{m=0}^n C_n^{(m)} = \emptyset$, then set $\tau \leftarrow n$, and declare a detection. 
        \end{itemize}
    (In the last step, we have implicitly used the nestedness discussed Remark~\ref{rem:nested}, but if the CSs are not nested, we can use the stopping criterion $\cap_{m=0}^n \cap_{i=m}^n C_i^{(m)} = \emptyset$.)
    \end{definition}
    Compared to the \bcsdetector of~\citet{shekhar23changepoint}, the main change in the above scheme is that it creates a new forward CS  with each new observation, instead of a new backward CS. We now show that our new SCD scheme defined above admits a nonasymptotic lower bound on the ARL. 
    \begin{theorem}[ARL control]
        \label{prop:arl-control}
        The changepoint detection scheme described in~\Cref{def:repeated-fcs-detector} controls the average run length~(ARL) at level $1/\alpha$. That is, when $T=\infty$, our proposed stopping time $\tau$ satisfies 
            $\mathbb{E}_{\infty}[\tau] \geq 1/\alpha$. %
    \end{theorem}
    The proof of this result is in~\Cref{proof:arl-control}. 
    Note that for the \bcsdetector, used with level-$(1-\alpha)$ CSs, \citet{shekhar23changepoint} obtained a lower bound of $1/2\alpha -3/2$ on the ARL. Thus, our \rfcsdetector achieves an improved (approximately by a factor of $2$) lower bound, while matching the detection delay guarantees of~\bcsdetector, as we show in the next section. 

    \begin{remark}
        \label{remark:PFA} 
        An alternative performance measure to ARL is the \emph{probability of false alarms~(PFA)}, which is  equal to the probability that the stopping time $\tau$ is finite; that is $\mathbb{P}_{\infty}(\tau<\infty)$. If we modify our \rfcsdetector to use a level-$(1- 6\alpha/(n^2 \pi^2))$ CS in rounds $n \geq 1$, then the resulting scheme ensures
        \begin{align}
            \mathbb{P}_{\infty}(\tau<\infty) \leq \sum_{n \geq 1} \mathbb{P}_{\infty} \lp\left\{ (C_t^{(n)})_{t \geq n} \text{ miscovers } \theta_0 \right\}\rp \leq  \alpha\sum_{n \geq 1}\frac{6}{\pi^2 n^2} = \alpha.
        \end{align}
        This implies that the ARL of the above modified \rfcsdetector scheme is infinity, since $\mathbb{E}_{\infty}[\tau] \geq (1-\alpha) \times \infty = \infty$. This significantly stronger control over false alarms comes at the cost of an increase in detection delay. In particular, we can show that for most CSs,  the detection delay of this modified scheme will have a logarithmic dependence on $T$. This means that the worst case~(over all $T$ values) detection delay of this scheme is usually unbounded. 
    \end{remark}

    \subsection{Analyzing the detection delay}
    \label{subsec:delay-analysis}

        We now state an assumption under which we will analyze the detection delay of our SCD scheme. 
        \begin{assumption}
            \label{assump:cs-width} Letting $d$ denote a metric on  $\Theta$, $X^n$ be shorthand for $(X_1,\dots,X_n)$, and $(C_n)_{n \geq 0}$ be a given confidence sequence, we assume that the width of the set $C_n \equiv C_n(X^n, \alpha)$ has a deterministic  bound
            \begin{align}
                 \sup_{\theta', \theta'' \in C_n} d(\theta', \theta'') \stackrel{a.s.}{\leq} w(n, P, \alpha), 
            \end{align}
            such that $\lim_{n \to \infty} w(n, P, \alpha) = 0$, for all  $P \in \mc{P}$ and $\alpha \in (0, 1]$. 
        \end{assumption}
            The above assumption requires the existence of a deterministic envelope function for the diameter of the confidence sequence, which converges to zero pointwise for every $(P, \alpha)$, as $n$ increases. This is a very weak assumption, and essentially all known CSs satisfy it. 
    
        We now analyze the detection delay  of our SCD scheme for~\Cref{problem:scd} under~\Cref{assump:cs-width}. 
        \begin{theorem}
            \label{prop:detection-delay-1}  Consider the SCD problem with observations $X_1, X_2, \ldots$ such that $(X_n)_{n \leq T}$ are drawn from a distribution $P_0$~(with parameter $\theta_0$), while $(X_n)_{n >T}$ are drawn from a product distribution $P_1$~(with parameter $\theta_1 \neq \theta_0$), and are independent of the pre-change observations. Suppose the \rfcsdetector from~\Cref{def:repeated-fcs-detector} is applied to this problem, with the CSs  satisfying~\Cref{assump:cs-width}. Let $\mc{E} := \{\theta_0 \in \cap_{n=1}^T C^{(1)}_n\}$ denote the ``good event'' (having at least $(1-\alpha)$ probability) that the first CS  covers the true parameter up to the changepoint. For~\Cref{problem:scd}~(non-partitioned), if  $T$ is large enough to ensure that $w(T, P_0, \alpha) < d(\theta_0, \theta_1)$, then the detection delay of our proposed scheme satisfies 
            \begin{align}
                &\mathbb{E}_T[(\tau-T)^+|\mc{E}] \leq \frac{3}{1-\alpha} u(\theta_0, \theta_1, T), \\
                \text{where} \quad 
                &u(\theta_0, \theta_1, T) \defined \min \{n \geq 1: w(T, P_0, \alpha) + w(n, P_1, \alpha) < d(\theta_0, \theta_1)\}. 
            \end{align}
            For~\Cref{problem:scd-separated}~(partitioned), with a known  $\Theta_0$, the \rfcsdetector satisfies 
            \begin{align}
                &\mathbb{E}_T[(\tau-T)^+]  \leq \frac{3}{1-\alpha} u(\Theta_0, \theta_1), \\
                \text{where} \quad 
                &u(\Theta_0, \theta_1) \defined \min \{n \geq 1:  w(n, P_1, \alpha) < \inf_{\theta' \in \Theta_0} d(\theta', \theta_1)\}. \label{eq:u-def-problem-2}
            \end{align}
            for all values of $T<\infty$. 
        \end{theorem}
        The proof of this result adapts the arguments developed by~\citet{shekhar23changepoint} for analyzing the \bcsdetector, and we present the details in~\Cref{proof:detection-delay-1}.         
        \begin{remark}
            \label{remark:arl-comparison}
            The above detection delay bound \emph{exactly} matches that obtained by the \bcsdetector of ~\citet{shekhar23changepoint}, which had an ARL guarantee of $\text{ARL} \geq 1/(2\alpha) - 3/2$. 
            Recalling~\Cref{prop:arl-control}, our new scheme achieves a significantly improved (by a factor of $2$) nonasymptotic lower bound on the ARL as compared to the \bcsdetector, while matching their detection delay performance. 
        \end{remark}
    
        This result provides a general detection delay bound applicable to a large class of problems. However, this guarantee can be sharpened when the \rfcsdetector is applied to  problems with additional structure. We demonstrate this next, through an application of our scheme to the problem of detecting changes in mean of bounded random variables.

    \subsection{A nonparametric example: change in mean for bounded random variables} 
    \label{subsec:bounded-mean}
        We now analyze the performance of our changepoint detection scheme the problem of detecting changes in the mean of bounded real-valued random variables supported on $\mc{X}=[0, 1]$. Note that despite the simple observation space, the class of distributions on this $\mc{X}$ is highly composite and nonparametric. In particular, there does not exist a common dominating measure for all distributions in this class, which renders likelihood based techniques inapplicable to this problem.

        Formally, we consider the instances of~\Cref{problem:scd} and~\Cref{problem:scd-separated} with $\mc{X} = [0,1]$, and the parameter space $\Theta = [0,1]$ with metric $d(\theta, \theta') = |\theta - \theta'|$. For~\Cref{problem:scd-separated}, we assume that the pre-change mean $\theta_0$ lies in a known set $\Theta_0 \subset \Theta$. For an unknown value $T \in \mathbb{N} \cup \{\infty\}$, the distribution generating the observations changes from $P_0$, with mean $\theta_0$, to another distribution $P_1$, with mean $\theta_1 \neq \theta_0$. The (unknown) change magnitude is denoted by $d(\theta_0, \theta_1) = \Delta = |\theta_1 - \theta_0|$. 
        
        For this problem, we employ our \rfcsdetector strategy using an instance of the betting-based construction of CSs for the means of bounded random variables~(details in~\Cref{proof:bounded-mean}) proposed by~\citet{waudby2023estimating}. Our next result analyzes its performance. 
        \begin{proposition}
            \label{prop:bounded-mean}
            \sloppy Consider the above problem of detecting changes in mean with bounded observations, under these additional conditions: (i) the post-change observations are independent of the pre-change observations, and (ii) both, the pre- and post-change observations are \iid~(that is $P_0, P_1$ are infinite products of some distributions $p_0, p_1$ on $\mc{X}$).  For~\Cref{problem:scd}~(non-partitioned), if $T \geq 64 \log(64/\Delta^2\alpha)/\Delta^2$,  the \rfcsdetector instantiated with the betting CS~(details in~\Cref{proof:bounded-mean}) satisfies: 
            \begin{align}
                &\mathbb{E}_{\infty}[\tau] \geq \frac{1}{\alpha}, \quad \text{and} \quad \mathbb{E}_{T}[(\tau - T)^+| \mc{E}] = \mc{O}\lp  \frac{\log(1/\alpha K_1)}{K_1}\rp, \label{eq:bounded-mean-1}\\
                \text{where} \quad &K_1 = K_1(P_1, \theta_0) \defined \inf_{P_{\theta}: |\theta - \theta_0| \leq \Delta/2} \dkl(p_1 \parallel p_{\theta}). 
            \end{align}
            In the display above, $\mc{E}$ is the ``good event'' in~\Cref{prop:detection-delay-1}, having probability at least $1-\alpha$. 

            For~\Cref{problem:scd-separated}~(partitioned),  the \rfcsdetector  satisfies the following:  
            \begin{align}
                &\mathbb{E}_{\infty}[\tau] \geq \frac{1}{\alpha}, \quad \text{and} \quad 
                \mathbb{E}_T[(\tau-T)^+] = \mc{O}\lp \frac{\log(1/\alpha K_2 )}{K_2} \rp, \label{eq:bounded-mean-2} \\ \text{where} \quad &K_2 \equiv K_2(P_1, \Theta_0) = \inf_{P_{\theta}: \theta \in \Theta_0} \dkl(p_1 \parallel p_{\theta}). 
            \end{align}
            In the statements above, $P_{\theta} = p_\theta^{\infty}$ denotes any product distribution on $\mc{X}^{\infty}$ with mean $\theta$. 
        \end{proposition}
        The proof of this result is in~\Cref{proof:bounded-mean}, and relies  on a careful analysis of the behavior of betting CSs. If the pre-change distribution $P_0$ is also \iid~(say $P_0 = p_0^{\infty}$), and is known, then $K_2$ in~\eqref{eq:bounded-mean-2} reduces to $\dkl(p_1 \parallel p_0)$. The resulting detection delay is order optimal, according to~\citet{lorden1971procedures}[Theorem 3], and furthermore, this optimality is achieved for an unknown $P_1$ lying in a nonparametric distribution class. 
        \begin{remark}
            By an application of Pinsker's inequality~\citep[Lemma 2.5]{10.5555/1522486}, we know that both $K_1$ and $K_2$ are $\Omega(\Delta^2)$, which gives us the weaker upper bound on the detection delay, $\mc{O}\lp \log(1/\alpha\Delta)/ \Delta^2 \rp$. This is the upper bound on the detection delay derived by~\citet{shekhar23changepoint} for the change of mean detection problem, using the empirical-Bernstein CS of~\citet{waudby2023estimating}, and a direct application of the general delay bound of~\Cref{prop:detection-delay-1}.
        \end{remark}

\section{Connection to Lorden's reduction from SCD to testing}
\label{sec:lorden} 
    Using the duality between confidence sequences and sequential hypothesis tests, we now show that our \rfcsdetector strategy is a generalization of a well-known result of~\citet{lorden1971procedures}, that reduces the problem of SCD (with separated  distribution classes) to that of repeated sequential tests. 
    
    Lorden's work built upon the interpretation of \cusum algorithm as repeated sequential probability ratio tests~(SPRT) for known pre- and post-change distributions by \citet{page1954continuous}. 
    In particular, \citet{lorden1971procedures} considered a parametric SCD problem with a known pre-change distribution $P_0$, and a parametric composite class of post-change distributions~$\{P_{\theta_1}: \theta_1 \in \Theta_1\}$. Then, given a  sequential test, or equivalently, extended stopping time, $\{N(\alpha): \alpha \in (0,1)\}$, satisfying $\mathbb{E}_{P_0}(N(\alpha) < \infty) \leq \alpha$, \citet{lorden1971procedures} proposed the following SCD strategy: 
    \begin{itemize}
        \item For every $m \geq 1$, define $N^{(m)}(\alpha)$ as the stopping rule $N(\alpha)$ applied to the observations $X_m, X_{m+1}, \ldots$. 
        \item Using these, declare the changepoint at the time $N^*(\alpha)$, defined as 
        \begin{align}
            N^*(\alpha) = \inf_{m \geq 1}\; \{ N^{(m)}(\alpha) + m - 1 \}. \label{eq:lorden-1}
        \end{align}
    \end{itemize}
    In words, this scheme can be summarized as: \emph{initiate a new sequential level-$\alpha$ test with every new observation, and stop and declare a detection as soon as one of the active tests rejects the null}. For this scheme, \citet{lorden1971procedures} established the ARL control; that is, $\mathbb{E}_{P_0}[N^*(\alpha)] \geq 1/\alpha$, for the specified $\alpha \in (0,1)$. Furthermore, under certain assumptions on the expected stopping time of the test $N(\alpha)$ under the alternative, \citet{lorden1971procedures} also established the minimax optimality of the scheme in the regime of $\alpha \to 0$.    

    Due to the duality between confidence sequences and sequential tests, it is easy to verify that the SCD scheme defined in~\eqref{eq:lorden-1} is a special instance of our general strategy. In particular, let us consider the problem studied by~\citet{lorden1971procedures}, with a known pre-change distribution $P_0$ with parameter $\theta_0$~(which could, for instance, be the distribution function of $P_0$), and a class of post-change distributions $\{P_{\theta_1}: \theta_1 \in \Theta_1\}$. Given a method for constructing  confidence sequences for $\theta$, denoted by $\mc{C}$,  we can define a sequential test 
    \begin{align}
        N(\alpha) = \inf \{n \geq 1: \theta_0 \not \in C_n\}, \quad \text{where} \quad C_n = \mc{C}(X_1, \ldots, X_n; \alpha). 
    \end{align}
    \sloppy By the uniform coverage guarantee of confidence sequences, we have $\mathbb{P}_{P_0}\lp N(\alpha)< \infty \rp = \mathbb{P}_{P_0} \lp \exists n \in \mathbb{N}: \theta_0 \not \in C_n \rp \leq \alpha$. Thus, $N(\alpha)$ is a valid level-$\alpha$ sequential test for the simple null $\{P_0\}$. 
    Similarly, $N^{(m)}$ for $m \geq 1$, can be defined as  the stopping time $N(\alpha)$ constructed using observations $(X_n)_{n \geq m}$ starting at time $m$. More specifically, we have 
    \begin{align}
        N^{(m)}(\alpha) = \inf \{n -m+1: \theta_0 \not \in C^{(m)}_n, \; n\geq m\}, \quad \text{where} \quad C_n^{(m)} = \mc{C}\lp X_m, X_{m+1} \ldots, X_n; \alpha \rp. 
    \end{align}
    Using this, we can define $N^*(\alpha) = \inf{m \geq 1} \{ N^{(m)}(\alpha) + m -1\}$. We now observe, that this $N^*(\alpha)$ is exactly the same as our proposed stopping time, $\tau$, in~\Cref{def:repeated-fcs-detector}, with $\Theta_0 = \{\theta_0\}$. In particular, for all $n \geq 1$, we have 
    \begin{align}
        \{N^*(\alpha) \leq n\} &= \{ \exists n' \leq n, \exists m \leq n': N^{(m)}(\alpha) = n'-m+1\} \\
        &= \{ \exists n' \leq n, \exists m \leq n': \theta_0 \not \in C^{(m)}_{n'} \}  
        = \{ \exists n' \leq n: \big( \cap_{m=1}^{n'} C^{(m)}_{n'} \big) \cap \{\theta_0\} = \emptyset \}  \\
        &= \{ \cap_{n'=1}^n \cap_{m=0}^{n'} C^{(m)}_{n'} = \emptyset \}  
         = \{ \cap_{m=0}^n \cap_{n'=m}^{n} C^{(m)}_{n'} = \emptyset \} = \{\tau \leq n\}. 
    \end{align}
    In the last two equalities, we have used the fact that $C_n^{(0)}$ for all $n \geq 0$ is equal to $\Theta_0 = \{\theta_0\}$. 
    Thus, our general strategy proposed in~\Cref{def:repeated-fcs-detector} subsumes the SCD scheme of~\citet{lorden1971procedures}. We end our discussion with the following two remarks: 
    \begin{itemize}
        \item While we instantiated our scheme above for the case of a singleton null, $\{P_0\}$, the exact same construction also applies to the case of a composite null $\{P_{\theta_0}: \theta_0 \in \Theta_0\}$, with $\Theta_0 \cap \Theta_1 = \emptyset$. The only modification needed is to update the stopping time $N(\alpha)$  to be equal to $\inf \{n \geq 1: \Theta_0 \cap C_n = \emptyset \}$. By~\Cref{prop:arl-control}, the resulting SCD scheme still controls the ARL at the required level $1/\alpha$. 

        \item Note that the e-detector framework, developed by~\citet{shin2022detectors}, also generalizes Lorden's scheme to work for composite, and nonparametric pre- and post-change distribution classes~($\mc{P}_0$ and $\mc{P}_1$ respectively). However, the e-detectors were developed explicitly for~\Cref{problem:scd-separated}~(partitioned); that is for a known class of pre-change distributions $\mc{P}_0$~(although the general idea could be suitably adapted for the non-separated formulation in some cases). This is unlike our scheme that is easily applicable to both the partitioned and non-partitioned formulations of the SCD problem. 
    \end{itemize}

\section{Deferred proofs from Section~\ref{sec:proposed-scheme}}
\label{sec:proofs} 
    In this section, we present the proofs of the three technical results stated in~\Cref{sec:proposed-scheme}. 
    
    \subsection{Proof of Theorem~\ref{prop:arl-control}}
    \label{proof:arl-control}
         We prove this statement in three steps. First, we define an e-process~(recalled in~\Cref{def:e-process}) corresponding to every confidence sequence $(C_n^{(m)})_{n \geq m}$ involved in our scheme.  Then, using these e-processes we introduce an e-detector $(M_n)_{n \geq 1}$, that is, a process adapted to the natural filtration $\filtration$ that satisfies $\mathbb{E}_{\infty}[M_{\tau'}] \leq \mathbb{E}_{\infty} [\tau']$ for all stopping times $\tau'$. Finally, we show that our stopping time $\tau$, introduced in~\Cref{def:repeated-fcs-detector}, is larger than $\tau' = \inf \{n \geq 1: M_n \geq 1/\alpha\}$, defined using the e-detector.   This allows us to leverage~\citet[Proposition~2.4]{shin2022detectors} to conclude that $\mathbb{E}_{\infty}[\tau'] \geq 1/\alpha$, which implies required statement about the ARL of $\tau$.

        Since we prove this result by attaching an e-process to every CS, we recall their definition below.
        \begin{definition}[e-processes]
            \label{def:e-process} Given a  class of probability measures $\mc{P}$, and a filtration $\filtration \equiv (\mc{F}_n)_{n \geq 1}$ defined on some measurable space, a $\mc{P}$-e-process is a collection of nonnegative random variables $(E_n)_{n \geq 1}$ adapted to $\filtration$, satisfying $\mathbb{E}_P[E_{\tau'}] \leq 1$ for all $P \in \mc{P}$, and for all stopping times $\tau'$~(adapted to the same filtration). 
        \end{definition}

        \paragraph{Step 1. Construct an e-process for every CS.} For every CS starting with the $m^{th}$ observation, denoted by $(C_n^{(m)})_{n \geq m}$, we associate a process defined as 
        \begin{align}
            E_n^{(m)} = 
                \begin{cases}
                    0, & \text{if } n < m, \text{ OR } \text{if } n \geq  m, \text{ and } \theta_0 \in C_n^{(m)}, \\
                    \frac{1}{\alpha}, & \text{if } n \geq  m, \text{ and } \theta_0 \not \in C_n^{(m)}. 
                \end{cases}
        \end{align}
        It is easy to verify that for every $m \geq 1$, the process $\{E_n^{(m)}: n \geq 1\}$ is an e-process:
        \begin{itemize}
            \item For every $n \geq 1$, the value of $E_n^{(m)}$ is $\mc{F}_n = \sigma(X_1, \ldots, X_n)$ measurable. 
            \item For any stopping time $\tau'$, adapted to the filtration $\filtration$, we have 
            \begin{align}
                \mathbb{E}_{\infty}[E_{\tau'}^{(m)}] &= \mathbb{E}_{\infty} \lb 0 \times \bone_{\tau' < m} + \frac{1}{\alpha} \times \bone_{\tau'\geq m} \bone_{\theta_0 \not \in C_{\tau'}^{(m)}} \rb  
                 = \frac{1}{\alpha} \times \mathbb{E}_{\infty} \lb  \bone_{\tau'\geq m} \bone_{\theta_0 \not \in C_{\tau'}^{(m)}} \rb  \\
                & \leq \frac{1}{\alpha} \times \mathbb{E}_{\infty} \lb  \bone_{\theta_0 \not \in C_{\tau'}^{(m)}} \rb 
                = \frac{1}{\alpha} \times \mathbb{P}_{\infty} \lp \theta_0 \not \in C_{\tau'}^{(m)} \rp 
                 \leq \frac{1}{\alpha} \times \alpha = 1. 
            \end{align}
        \end{itemize}
        The last inequality uses the fact that $(C_n^{(m)})_{m \geq n}$ is a level-$(1-\alpha)$ CS for $\theta_0$. Thus, for every $m \geq 1$, the process $(E_n^{(m)})_{n \geq 1}$  is a valid e-process.  

        \paragraph{Step 2. Construct an e-detector.} For every $n \geq 1$, we define $M_n$ to be equal to $\sum_{m=1}^n E_n^{(m)}$, and observe that the process $(M_n)_{n \geq 1}$ is an \emph{e-detector}, as defined by~\citet[Definition 2.2]{shin2022detectors} because it satisfies the following two properties: 
        \begin{itemize}
            \item $(M_n)_{n \geq 1}$ is adapted to $(\mc{F}_n)_{n \geq 1}$: since for any $n \geq 1$, all the $E_n^{(m)}$ are $\mc{F}_n$-measurable by construction. 
            \item For any stopping time $\tau'$, we have 
                $\mathbb{E}_{\infty}[M_{\tau'}] \leq \mathbb{E}_{\infty}[\tau']$, as noted by~\citet[Definition~2.6]{shin2022detectors}.  
        \end{itemize}
       
        \paragraph{Step 3. Bound the ARL using the e-detector.} Finally, we translate the stopping criterion of our proposed scheme (stated as the non-intersection of the confidence sequences) in terms of the e-detector $(M_n)_{n \geq 1}$. In particular, we have 
        \begin{align}
            \{\tau \leq n\} = \{ \cap_{m=0}^n C_n^{(m)} = \emptyset \} 
            \subset \{ \exists m \in [n]: \theta_0 \not \in C_n^{(m)} \} 
            = \{ \exists m \in [n]: E_n^{(m)} = 1/\alpha \} 
            = \{M_n \geq 1/\alpha\}.  \label{eq:e-detector-1}
        \end{align}
        In words, when $T=\infty$, if the intersection of the CSs is empty prior to some time $n$, it means that at least of of the CSs constructed prior to $n$ must miscover. This in turn implies that the value of at least one of the e-processes at $n$ is equal to $1/\alpha$; or the value of the e-detector $M_n$ is at least $1/\alpha$. 
        Recall, that in the first equality above, we have assumed that the sets in the confidence sequences are nested; that is, $C_n^{(m)} \subset C_{n'}^{(m)}$ for every $m \leq n' < n$. This allows us to look only at the intersection of the most recent sets to define the stopping condition. 
        We now define a new stopping time 
        \begin{align}
            \tau' = \inf \{n \geq 1: M_n \geq 1/\alpha \}, 
        \end{align}
        and observe that it is stochastically dominated by $\tau$; that is,~\eqref{eq:e-detector-1} implies that 
        \begin{align}
            \{\tau' > n\} = \{M_n < 1/\alpha\} \subset \{\tau > n\}, \quad \text{which implies} \quad \mathbb{E}_{\infty}[\tau'] \leq \mathbb{E}_{\infty}[\tau]. 
        \end{align}
        From~\citet[Proposition~2.4]{shin2022detectors}, we know that $\mathbb{E}_{\infty}[\tau'] \geq 1/\alpha$, and we conclude the result by noting that $\mathbb{E}_{\infty}[\tau]  
 \geq \mathbb{E}_{\infty}[\tau']$ since $\tau$ stochastically dominates $\tau'$. \hfill \qedsymbol

    \subsection{Proof of Theorem~\ref{prop:detection-delay-1}}
    \label{proof:detection-delay-1}
        The proof of this result follows the general argument developed by~\citet{shekhar23changepoint} for analyzing their \bcsdetector strategy, with some modifications due to the use of forward CSs~(instead of backward CSs used in the \bcsdetector).  
    
        In particular, we consider blocks of the post-change observations, each of length $u \equiv u(\theta_0, \theta_1, T)$, starting at time $T_j = T + ju$ for $j \geq 0$. Note that all these blocks are independent of each other~(since $P_1$ is a product distribution), and also independent of the event $\mc{E} = \{ \theta_0 \in C^{(1)}_T\}$. Now, observe that for $k \geq 1$, we have $\{\tau > T_k \} = \cap_{j=1}^k \{ \tau > T_j \}$, which furthermore implies
        \begin{align}
            \{ \tau > T_k \} \cap \mc{E} \subset \cap_{j=1}^{k} \{ C_{T_j}^{(T_{j-1})} \cap C_T^{(1)} \neq \emptyset \} \cap \mc{E} \subset \cap_{j=1}^k \{C^{(T_{j-1})}_{T_{j}} \text{ miscovers } \theta_1 \} \cap \mc{E}. 
        \end{align}
        The last inclusion follows from the definition of $u \equiv u(\theta_0, \theta_1, T)$, and the event $\mc{E}$. Using this, we obtain:
        \begin{align}
            \sum_{t=T_k+1}^{T_{k+1}} \mathbb{P}_T \lp \tau \geq t | \mc{E} \rp  &\leq \sum_{t=T_k+1}^{T_{k+1}} \mathbb{P}_T \lp \tau \geq T_k| \mc{E} \rp = u \mathbb{P}_T\lp \tau > T_k | \mc{E} \rp  
             \leq u \times \mathbb{P}_T \lp  \cap_{j=1}^k \{C^{(T_{j-1})}_{T_{j}} \text{ miscovers } \theta_1 \} \cap \mc{E} |\mc{E} \rp \\
             & \stackrel{(i)}{\leq} \frac{u}{\mathbb{P}_T(\mc{E})} \times \prod_{j=1}^k \mathbb{P}_T \lp   \{C^{(T_{j-1})}_{T_{j}} \text{ miscovers } \theta_1 \} \cap \mc{E} \rp \leq \frac{u \alpha^k}{1-\alpha}. 
        \end{align}
        The inequality $(i)$ uses the fact that $\mc{E}$ only depends on the pre-change observations, and hence is independent of the post-change CSs. 
        We now fix some $k_0>1$, and observe that 
        \begin{align}
            \mathbb{E}_T[(\tau-T)+] \leq k_0  u  + \sum_{k=k_0}^{\infty} \frac{u \alpha^k}{1-\alpha} = u \lp k_0 + \frac{\alpha^{k_0}}{1-\alpha} \rp. 
        \end{align}
        By setting $k_0 = \lceil \log(1/1-\alpha)/\log(1/\alpha) \rceil$, we get the required statement for~\Cref{problem:scd}.

        To prove the second part of~\Cref{prop:detection-delay-1}, we proceed as above, considering blocks of post-change observations of length $u \equiv u(\Theta_0, \theta_1)$ as defined in~\eqref{eq:u-def-problem-2}. We then define $T_j = T + ju$ for $j \geq 0$, and note that 
        \begin{align}
            \{\tau > T_k\} \subset \cap_{j=1}^k \{ C^{(T_{j-1})}_{T_j} \cap \Theta_0 \neq \emptyset \}  
            \subset \cap_{j=1}^k \{C^{(T_{j-1})}_{T_{j}} \text{ miscovers } \theta_1 \}. 
        \end{align}
        The rest of the argument, then proceeds exactly as in the first part. 

    \subsection{Proof of Proposition~\ref{prop:bounded-mean}}
    \label{proof:bounded-mean}
        Before presenting the proof of~\Cref{prop:bounded-mean}, we first recall some of details of the betting CS first proposed by~\citet{waudby2023estimating}. 
        
        \paragraph{Background on betting CS.} Given observations $X_1, X_2, \ldots$ drawn from an independent process with mean $\theta$, the betting CS is defined as 
        \begin{align}
            C_n = \{s \in [0,1]: W_n(s) < 1/\alpha\}, \quad \text{with} \quad W_n(s) \defined \prod_{i=1}^n (1 + \lambda_t(s) (X_t-s)), \quad \text{for all } s \in [0,1], 
        \end{align}
        where $\{\lambda_t(s): t \geq 1, s \in [0,1]\}$ are predictable bets, taking values in $[-1/(1-s), 1/s]$. For certain betting strategies, such as the \emph{mixture method}~\citep[\S~4.3]{hazan2016introduction}, the \emph{regret} is logarithmic for all $s$. In particular, this implies that 
        \begin{align}
             \sup_{ \lambda \in \lb \frac{-1}{1-s}, \frac{1}{s}\rb }\sum_{t=1}^n \log \lp 1 + \lambda (X_t - s) \rp  - \sum_{t=1}^n \log\lp 1 + \lambda_t(s)(X_t-s) \rp \leq 2 \log n, \quad \text{for all } n \geq 13.
        \end{align}
        Note that this idea of using the mixture method with known regret guarantees, for the specific context of betting CS was first considered by~\citet{orabona2021tight}. 
        We now present the details of the proof of~\Cref{prop:bounded-mean}. First we show that under the condition that $T \geq 64\log(64/\Delta^2\alpha)/\Delta^2$, the analysis of the first setting~(i.e., with unknown $\Theta_0$) can be reduced to the second case~(with known $\Theta_0$). Then, we present the details of the proof of the second setting. 

        \paragraph{Proof of~(\ref{eq:bounded-mean-1}).}
            Using the fact that $\log(1+x) \geq x - x^2/2$, we can further lower bound 
            $\log W_n(s)$ with 
            \begin{align}
                \log W_n(s) & \geq \sup_{ \lambda \in \lb \frac{-1}{1-s}, \frac{1}{s}\rb }\sum_{t=1}^n \lambda (X_t-s) - \frac{\lambda^2}{2}(X_t-s)^2  - \log(n^2). 
            \end{align}
            By setting the value of $\lambda$ to $\frac{1}{n}\sum_{t=1}^n X_t-m$, and on simplifying, we can show that the betting CS after $n$ observation satisfies $|C^{(1)}_n| \leq 4 \sqrt{\log(n/\alpha) /n}$.  This implies that for $T \geq 64 \log(64/\Delta^2 \alpha)/\Delta^2$, the width of the CS starting at time $1$ must be smaller than $\Delta/2 = |\theta_1-\theta_0|/2$. If the event $\mc{E} = \{ \theta_0 \in C^{(1)}_T \}$ happens~(recall that this is a probability $1-\alpha$ event), then we know that $\theta_0 \in \widetilde{\Theta}_0 \defined \{\theta: |\theta-\theta_0|\leq \Delta/2\}$. This set $\widetilde{\Theta}_0$ plays the role of the known pre-change parameter class in the analysis. Hence the rest of the proof to obtain the upper bound stated in~\eqref{eq:bounded-mean-1} proceeds exactly as in the case when the pre-change distribution class is known, and we present the details for the latter case next.

        \paragraph{Proof of~(\ref{eq:bounded-mean-2}).} Since the proof of this result is long, we break it down into four simpler steps. 
        \bigskip
    
            \emph{Step 1: Bound $(\tau - T)$ with the maximum of a class of stopping times $(N_{\theta})_{\theta \in \Theta_0}$.}
            Introduce the stopping times $N_m = \inf \{n \geq m: C^{(m)}_n \cap \Theta_0 = \emptyset \}$, and note that 
            \begin{align}
                \tau \leq \min_{m \geq 1} N_m, \quad \text{which implies that} \quad  \mathbb{E}_T[\tau] \leq \min_{m \geq T+1} \mathbb{E}_T[N_m].
            \end{align}
            Since the post change observations are assumed to be drawn \iid from a distribution with mean $\theta_1$, all the $N_m$ for $m \geq T+1$, have the same distribution, and thus, the same expected value. Hence, it suffices to get an upper bound on $N_{T+1}$. 
    
            To simplify the ensuing argument, we will use $C_n$ to denote $C_{n+T+1}^{(T+1)}$ for $n \geq 1$. Furthermore, we also assume that $\theta_1 < \theta$ for all $\theta \in \Theta_0$, and $\inf_{\theta \in \Theta_0} \theta = \theta_1 + \Delta$.  
            By the definition of betting CS, the stopping time $N_{T+1}$ can be written as the supremum of a collection of stopping times: 
            \begin{align}
                N_{T+1} - T = \sup_{\theta \in \Theta_0} N_\theta, \quad \text{where} \quad N_{\theta} = \inf \{n \geq 1: W_n(\theta) \geq 1/\alpha\}. 
            \end{align}
    
            \emph{Step 2: Bound $(N_\theta)_{\theta \in \Theta_0}$ with a monotonic class of stopping times.}
            Next, we will upper bound each $N_{\theta}$ with another stopping time $\gamma_\theta$, which have the property that $\gamma_{\theta'} < \gamma_{\theta}$ for $\theta' > \theta$. In particular, using the regret guarantee of the betting strategy, observe the following: 
            \begin{align}
                \log(W_n(\theta)) &\geq  \sup_{\lambda \in \lb \frac{-1}{1-\theta}, \frac{1}{\theta} \rb } \sum_{t=1}^n \log(1 + \lambda(X_t-\theta)) - 2 \log n \\
                &\geq  \sup_{\lambda \in\blue{\lb 0, \frac{1}{1-\theta} \rb}} \sum_{t=1}^n \log(1 + \lambda\blue{(\theta - X_t)} - 2 \log n \defined Z_n(\theta) - 2\log n.  
            \end{align}
            Define a new stopping time $\gamma_\theta = \inf \{n \geq 1: Z_n(\theta) -2 \log n \geq \log(1/\alpha)\}$, and note that the above display implies $\gamma_\theta \geq N_\theta$, and thus we have $N_{T+1}-T \leq \sup_{\theta \in \Theta_0} \gamma_{\theta}$. We now show the monotonicity of $\gamma_\theta$. 
    
            For any $\theta' > \theta$, we have $ \lambda(\theta'-X_t) \geq \lambda (\theta - X_t)$ for any $\lambda >0$, which implies  that $\sum_{t=1}^n \log (1 + \lambda(\theta'-X_t)) \geq \sum_{t=1}^n \log(1 + \lambda(\theta-X_t))$. Thus, we have the following relation (for $\theta' \geq \theta$): 
            \begin{align}
                Z_n(\theta') &= \sup_{\lambda \in \lb 0, \frac{1}{1-\theta'} \rb } \sum_{t=1}^n \log( 1+ \lambda(\theta' - X_t) ) 
                \geq \sup_{\lambda \in \lb 0, \blue{\frac{1}{1-\theta}} \rb } \sum_{t=1}^n \log( 1+ \lambda(\theta' - X_t) ) \\
                &\geq \sup_{\lambda \in \lb 0, \frac{1}{1-\theta} \rb } \sum_{t=1}^n \log( 1+ \lambda(\blue{\theta} - X_t) )  = Z_n(\theta). 
            \end{align}
            Thus, $Z_n(\theta') \geq Z_n(\theta)$, which implies that $\gamma_{\theta'} \leq \gamma_{\theta}$, and in particular, $\gamma_{\theta} \leq \gamma_{\theta_1 + \Delta}$ for all $\theta \in \Theta_0$. This leads to the required conclusion 
            \begin{align}
                N_{T+1} - T \leq \sup_{\theta \in \Theta_0} N_{\theta} \leq \sup_{\theta \in \Theta_0} \gamma_{\theta} \leq \gamma_{\theta_1 + \Delta}. 
            \end{align}
            This is a crucial step, as it reduces the task of analyzing the supremum of a large collection of stopping times, into that of analyzing a single stopping time $\gamma_{\theta_1 + \Delta}$. 
    
            \emph{Step 3: Bound $\gamma_{\theta_1 + \Delta}$ with the `oracle' stopping time $\rho^*$.} Let $\lambda^* \equiv \lambda^*(\theta_1 + \Delta)$ denote the log-optimal betting fraction, defined as $\argmax_{\lambda \in [0, 1/(1-\theta_1-\Delta)} \mathbb{E}[\log(1+\lambda(\theta_1 + \Delta - X))]$, where $X$ is drawn from the post-change distribution. By definition then, we have 
            \begin{align}
                Z_n(\theta_1 + \Delta) \geq Z_n^*(\theta_1 + \Delta) \defined \sum_{t=1}^n \log(1 + \lambda^*(\theta_1 + \Delta - X_t)), 
            \end{align}
            which immediately implies 
            \begin{align}
                \gamma_{\theta_1 + \Delta} \leq \rho^* \defined \inf \{n \geq 1: Z_n^*(\theta_1 + \Delta) \geq \log(n^2/\alpha) \}. 
            \end{align}
            The stopping time $\rho^*$ is much easier to analyze as it is the first crossing of the boundary $\log(n^2/\alpha)$ by the random walk $Z_n^*(\theta_1 + \Delta)$ with \iid increments. 
    
            \emph{Step 4: Evaluate the expectation of $\rho^*$.} Observe that $Z_n^* \equiv Z_n^*(\theta_1 + \Delta) = \sum_{t=1}^n V_t$, with $V_t = \log(1+ \lambda^*(\theta_1 + \Delta - X_t))$. Without loss of generality, we can assume that $\lambda^* < 1/(1-\theta_1 - \Delta)$~(if not, we simply repeat the argument with $\lambda^*-\epsilon$ for an arbitrarily small $\epsilon>0$), and hence $(V_t)_{t\geq 1}$ are \iid and bounded increments, which means that $\mathbb{E}[V_t]<\infty$. In fact, by the dual definition of the information projections~\citep{honda2010asymptotically}, we have $\mathbb{E}[V_t] = K_2 \equiv K_2(P_1, \Theta_0)$. 
            Next, with $n_0 \defined \inf \{n \geq 1: \log(n^2/\alpha)/n < K_2/2\}$, we have  for $n \geq n_0$ by an application of Hoeffding's inequality: 
            \begin{align}
                \mathbb{P}\lp \rho^* > n \rp \leq \mathbb{P}\lp \frac{1}{n} \sum_{t=1}^n V_t - K_2 \leq - \frac{K_2}{2} \rp \leq \exp \lp -c'' n \rp, 
            \end{align}
            for some $c''>0$. Hence, the expectation of $\rho^*$ satisfies 
            \begin{align}
                \mathbb{E}[\rho^*] = \sum_{n \geq 0} \mathbb{P}\lp \rho^* > n \rp \leq n_0 + \sum_{n \geq n_0} \exp \lp - c'' n \rp = n_0 + \frac{e^{-c'' n_0}}{1 - e^{-c''}} < \infty. 
            \end{align}
            Thus, both $\rho^*$ and $(V_t)_{t\geq 1}$ have bounded expectations, and we can appeal to Wald's lemma~\citep[Theorem~2.6.2]{durrett2019probability} to obtain $\mathbb{E}[Z^*_{\rho^*}] = \mathbb{E}[\rho^*] K_2$. Furthermore, by the definition of $\rho^*$, and the boundedness of $(V_t)_{t \geq 1}$, we can upper bound $\mathbb{E}[Z_{\rho^*}^*]$ with  $\log(1/\alpha) + 2\log(\mathbb{E}_T[\rho^*]) + c'$, where $c' = \max \{ \log(1+\lambda^*_{\theta_1 + \Delta}), \, \log(1-\lambda^*_{\theta_1 + \Delta})\}$. In other words, we have 
            \begin{align}
                \mathbb{E}[\rho^*] \leq \frac{ \log(1/\alpha) + 2\log(\mathbb{E}[\rho^*]) + c'}{K_2},  \quad \text{which implies} \quad 
                \mathbb{E}[\rho^*] = \mc{O}\lp \frac{\log(1/\alpha K_2)}{K_2} \rp. 
            \end{align}
            This completes the proof. 
 
\section{Conclusion}    
\label{sec:conclusion}
    In this paper, we proposed a changepoint detection scheme that constructs a new CS with every observation, and declares a detection as soon as the intersection of the active CSs becomes empty. The design of our scheme was motivated by the \bcsdetector of~\citet{shekhar23changepoint}, which proceeds by initializing new ``backward CSs'' with each new observation. We showed that our new scheme matches the detection delay performance of \bcsdetector, while improving the ARL lower bound by a factor of $2$. Furthermore, our scheme achieves this improvement under weaker model assumptions~(i.e., without needing the ability to construct CSs in both forward and backward directions). Interestingly, our proposed scheme can be seen as a nonparametric generalization of Lorden's reduction from SCD to repeated sequential testing, due to the duality between sequential testing and CSs. 

\subsection*{Acknowledgement}
The authors acknowledge support from NSF grants IIS-2229881 and DMS-2310718
\newpage 
\bibliographystyle{abbrvnat}
\bibliography{ref}

\begin{thebibliography}{14}
\providecommand{\natexlab}[1]{#1}
\providecommand{\url}[1]{\texttt{#1}}
\expandafter\ifx\csname urlstyle\endcsname\relax
  \providecommand{\doi}[1]{doi: #1}\else
  \providecommand{\doi}{doi: \begingroup \urlstyle{rm}\Url}\fi

\bibitem[Durrett(2019)]{durrett2019probability}
R.~Durrett.
\newblock \emph{Probability: Theory and Examples}.
\newblock Cambridge Series in Statistical and Probabilistic Mathematics.
  Cambridge University Press, 5th edition, 2019.

\bibitem[Hazan(2016)]{hazan2016introduction}
E.~Hazan.
\newblock Introduction to online convex optimization.
\newblock \emph{Foundations and Trends in Optimization}, 2\penalty0
  (3-4):\penalty0 157--325, 2016.

\bibitem[Honda and Takemura(2010)]{honda2010asymptotically}
J.~Honda and A.~Takemura.
\newblock An asymptotically optimal bandit algorithm for bounded support
  models.
\newblock In \emph{The Twenty Third Annual Conference on Learning Theory},
  pages 67--79. PMLR, 2010.

\bibitem[Lorden(1971)]{lorden1971procedures}
G.~Lorden.
\newblock Procedures for reacting to a change in distribution.
\newblock \emph{The Annals of Mathematical Statistics}, pages 1897--1908, 1971.

\bibitem[Orabona and Jun(2024+)]{orabona2021tight}
F.~Orabona and K.-S. Jun.
\newblock Tight concentrations and confidence sequences from the regret of
  universal portfolio.
\newblock \emph{IEEE Transactions on Information Theory (to appear)}, 2024+.

\bibitem[Page(1954)]{page1954continuous}
E.~S. Page.
\newblock Continuous inspection schemes.
\newblock \emph{Biometrika}, 41\penalty0 (1/2):\penalty0 100--115, 1954.

\bibitem[Shekhar and Ramdas(2023)]{shekhar23changepoint}
S.~Shekhar and A.~Ramdas.
\newblock Sequential changepoint detection via backward confidence sequences.
\newblock In \emph{Proceedings of the 40th International Conference on Machine
  Learning}, volume 202 of \emph{Proceedings of Machine Learning Research},
  pages 30908--30930. PMLR, 23--29 Jul 2023.

\bibitem[Shewhart(1925)]{shewhart1925application}
W.~A. Shewhart.
\newblock The application of statistics as an aid in maintaining quality of a
  manufactured product.
\newblock \emph{Journal of the American Statistical Association}, 20\penalty0
  (152):\penalty0 546--548, 1925.

\bibitem[Shewhart(1930)]{shewhart1930economic}
W.~A. Shewhart.
\newblock Economic quality control of manufactured product.
\newblock \emph{Bell System Technical Journal}, 9\penalty0 (2):\penalty0
  364--389, 1930.

\bibitem[Shin et~al.(2024)Shin, Ramdas, and Rinaldo]{shin2022detectors}
J.~Shin, A.~Ramdas, and A.~Rinaldo.
\newblock E-detectors: a nonparametric framework for online changepoint
  detection.
\newblock \emph{New England Journal of Statistics and Data Science (to
  appear)}, 2024.

\bibitem[Shiryaev(1963)]{shiryaev1963optimum}
A.~N. Shiryaev.
\newblock On optimum methods in quickest detection problems.
\newblock \emph{Theory of Probability \& Its Applications}, 8\penalty0
  (1):\penalty0 22--46, 1963.

\bibitem[Tartakovsky et~al.(2014)Tartakovsky, Nikiforov, and
  Basseville]{tartakovsky2014sequential}
A.~Tartakovsky, I.~Nikiforov, and M.~Basseville.
\newblock \emph{{Sequential analysis: Hypothesis testing and changepoint
  detection}}.
\newblock CRC Press, 2014.

\bibitem[Tsybakov(2009)]{10.5555/1522486}
A.~B. Tsybakov.
\newblock \emph{Introduction to Nonparametric Estimation}.
\newblock Springer Series in Statistics, 2009.

\bibitem[Waudby-Smith and Ramdas(2023)]{waudby2023estimating}
I.~Waudby-Smith and A.~Ramdas.
\newblock Estimating means of bounded random variables by betting.
\newblock \emph{Journal of the Royal Statistical Society B}, 2023.

\end{thebibliography}

\end{document}